\documentclass{amsart}

\usepackage[english]{babel}

\usepackage{fullpage}
\usepackage{geometry}
\usepackage{xcolor}
\usepackage{graphicx}
\usepackage{multirow}
\usepackage{wrapfig}
\usepackage{longtable}
\usepackage{float}
\usepackage{setspace} 
\usepackage{array}
\usepackage{booktabs,tabularx}
\usepackage{amsmath,amssymb,enumerate,url,amsthm}
\usepackage{ytableau}
\usepackage[numbers,sort]{natbib}
\usepackage{subcaption}
\usepackage{gensymb}
\usepackage{enumitem}
\usepackage{pgf,varwidth}
\usepackage{chngcntr}
\usepackage{hyperref,cases}

\usepackage[bottom]{footmisc}

\definecolor{dblue}{RGB}{0,54,135}
\definecolor{lorange}{RGB}{255,155,0}

\newcommand{\cp}{\mathrm{cp}}
\newcommand{\CP}{\mathcal{CP}}

\newcommand{\N}{\mathbb{N}}

\theoremstyle{definition}
\newtheorem*{remark}{Remark}
\newtheorem{theorem}{Theorem}
\newtheorem{cor}[theorem]{Corollary}

\newtheorem*{defn}{Definition}
\newtheorem{example}[theorem]{Example}

\numberwithin{theorem}{section}

\allowdisplaybreaks

\title{Copartitions}
\author{Hannah E. Burson and Dennis Eichhorn}
\date{\today}

\begin{document}
\maketitle

\begin{abstract}
    We develop the theory of copartitions, which are a generalization of partitions with connections to many classical topics in partition theory, including Rogers-Ramanujan partitions, theta functions, mock theta functions, partitions with parts separated by parity, and crank statistics.  
    Using both analytic and combinatorial methods, we give two forms of the three-parameter generating function, and we study several special cases that demonstrate the potential broader impact the study of copartitions may have.   
\end{abstract}

\section{Introduction}
In a study of combinatorial interpretations of Watson's third order mock theta functions, Andrews \cite{Andrews18} introduced the function  $\mathcal{EO}^*(n)$, which counts the number of integer partitions of $n$ with all even parts smaller than all odd parts, where the only part appearing an odd number of times is the largest even part. In introducing this function, Andrews opened the door to a wide array of further questions about partitions with parts separated by parity, first formalized by Andrews in \cite{AndrewsParity}. 

In, \cite{Andrews18},
Andrews noted many interesting properties of $\mathcal{EO}^*(n)$, including that its generating function is simply  $\frac{1}{2}(\nu(q)+\nu(-q))$, where $$\nu(q):=\sum_{n=0}^\infty \frac{q^{n^2+n}}{(-q;q^2)_{n+1}}$$ is one of Watson's third order mock theta functions. Additionally, Andrews noted that, like the ordinary partition function, 
$\mathcal{EO}^*(n)$ enjoys divisibility by five in an arithmetic progression:
$$
\mathcal{EO}^*(10n+8) \equiv 0 \pmod 5.
$$

In \cite{Chern19}, Chern provided a combinatorial proof of the generating function for $\mathcal{EO}^*(n)$. Then, in \cite{Chern21}, Chern studied further properties of $\mathcal{EO}^*(n)$, with a special focus on the remainder of the largest even part modulo $4$.

We note that here
we have adopted Chern's notation, rewriting Andrews' original $\overline{\mathcal{EO}}(n)$ from \cite{Andrews18} as $\mathcal{EO}^*(n)$. Chern \cite{Chern19} renamed $\overline{\mathcal{EO}}(n)$ to more easily discuss its overpartition analogs.

In this paper, we generalize $\mathcal{EO}^*(n)$ by introducing new partition-theoretic objects called \emph{copartitions}. These objects reveal an inherent symmetry in partitions counted by $\mathcal{EO}^*(n)$ that was not previously obvious.
Copartitions are counted by the
function $\cp_{a,b,m}(n)$ defined in Section \ref{sec:copartitionsIntro}, where $\cp_{1,1,2}(n)=\mathcal{EO}^*(2n)$. 
 Additionally, copartitions are
 in bijective correspondence with the subset of partitions called capsids, as defined by Garvan and Schlosser \cite{Capsids} 
to combinatorially treat the relationship between Ramanujan's tau function and $t$-core partitions.

We find that $\cp_{a,b,m}(n)$ has a beautiful infinite-product generating function (Section \ref{sec:copartitionsIntro}). 
Special cases of $\cp_{a,b,m}(n)$ provide a new framework for understanding a surprising number of classical partition-theoretic objects, 
including the Rogers-Ramanujan functions, the summatory function of $p(n)$, the function summing the largest parts over all partitions of $n$, and the function summing the perimeters over all partitions of $n$ (Section \ref{sec:cop:cases}).

\section{Background on Partitions}\label{sec:partitions}
In this paper, we use the language of integer partitions and $q$-series.  

By a partition of a non-negative integer $n$, we refer to a multiset of positive integers $\{\lambda_1,\lambda_2,\ldots,\lambda_r\}$ such that $\lambda_1\ge \lambda_2\ge \ldots \ge \lambda_r$ and $\lambda_1+\lambda_2+\ldots+\lambda_r=n$. We call $n$ the \emph{size} of the partition and write $|\lambda|=n$. When a partition $\lambda$ has many parts, it is useful to write $\lambda=\{1^{f_1},2^{f_2},3^{f_3},\ldots\}$ where only finitely many $f_i$ are nonzero and $\lambda$ has $f_i$ parts of size $i$. This notation is called frequency notation. For example, the partition $\{6,6,6,5,5,4,3,3,3,1,1,1\}$ can be written in frequency notation as $\{1^3,3^3,4,5^2,6^3\}$. When studying integer partitions, it can be useful to use the $q$-rising factorial $$(a;q)_n=\prod_{k=0}^{n-1}(1-aq^k).
$$ When $|q|<1$, we can also write $$(a;q)_\infty=\lim_{n\to \infty}(a;q)_n=\prod_{k=0}^\infty (1-aq^k).$$

Integer partitions can be graphically represented as an array of boxes (also called cells), with each row left-justified and each column top-justified. Due to the justification of the columns, the rows are in non-increasing order, so we can think of $\lambda_i$ as the number of boxes in the $i$th row of the diagram. This graphical representation is called a Young diagram (or Ferrers diagram). For a partition $\lambda$, we define the conjugate partition $\lambda'$ as the partition obtained by reflecting the Young diagram of $\lambda$ about the line $y=-x$. Equivalently, the part $\lambda'_i$ is equal to the number of boxes in the $i$th column of the Young diagram of $\lambda$.

Sometimes, especially when studying partitions with parts from a specific arithmetic progression, it is useful to work with an adaptation of a Young diagram called an $m$-modular diagram. To obtain the $m$-modular diagram for a partition $\lambda$, we represent parts of size $mk+r$ as a box containing an $r$ followed by $k$ boxes, each containing an $m$. Note that the first column of an $m$-modular diagram shows the remainders modulo $m$ of each part of the partition.

A reoccurring theme throughout this manuscript is the importance of the number of distinct parts of a partition. Thus, we propose using the term \emph{diversity} of a partition to mean the number of different part sizes occurring in that partition. We use $\mathrm{dv}(\lambda)$ to denote the diversity of the partition $\lambda$. In \cite{Kim2012}, Kim constructs an elaborate generating function for $p(k,n)$, the number of partitions of $n$ with diversity $k$, as a sum over partitions of $k$.

\section{Definition, Diagrams, and Generating Function}\label{sec:copartitionsIntro}

In this section, we introduce $(a,b,m)$-copartitions and the function $\cp_{a,b,m}(n).$
\begin{defn}
An $(a,b,m)$-copartition is a triple of partitions $(\gamma,\rho,\sigma)$, where each of the parts of $\gamma$ is at least $a$ and congruent to $a\pmod m$, each of the parts of $\sigma$ is at least $b$ and congruent to $b\pmod m$, and $\rho$ has the same number of parts as $\sigma$, each of which have size equal to $m$ times the number of parts of $\gamma$.\\
 When $a,b,m\ge 1$, we let $\cp_{a,b,m}(n)$ denote the number of $(a,b,m)$-copartitions of size $n$. 
\end{defn}

Although $\rho$ is completely determined by $\gamma$ and $\sigma$, the graphical representation we use suggests that the natural way to write the triple is $(\gamma,\rho,\sigma)$.  
We refer to $\gamma$ as the \emph{ground} of a copartition and $\sigma$ as the \emph{sky}.

\begin{example}
The $(1,3,4)$-copartitions of size 12 are
   $$ \left (\{9,1^3\},\emptyset,\emptyset \right ),\;
   \left (\{5^2,1^2\},\emptyset,\emptyset \right ),\; \left (\{5,1^7\},\emptyset,\emptyset \right ),  \left (\{1^{12}\},\emptyset,\emptyset \right ),\;$$ 
   $$\left (\{5\},\{4\},\{3\} \right ),\; \left (\{1\},\{4\},\{7\} \right ), \; \text{and } \left (\emptyset,\emptyset, \{3^4\} \right ).
$$
Thus, $\cp_{1,3,4}(12)=7$. 
\end{example}
  
To represent the $(a,b,m)$ copartition $(\sigma, \rho, \gamma)$ graphically, we append the $m$-modular diagram for $\sigma$ to the right of the $m$-modular diagram for $\rho$. Then, we append the conjugate of the $m$-modular diagram for $\gamma$ below $\rho$. 

\begin{example} The following diagram represents the $(a,b,m)$-copartition $(\{3m+a,2m+a,2m+a,a\},\{4m,4m\},\{3m+b,2m+b\})$.

\begin{center}
\begin{minipage}{0.35\textwidth}
\beginpgfgraphicnamed{copartitions-f1}
\begin{tikzpicture}[remember picture]

\node (n) {\begin{varwidth}{5cm}{
\begin{ytableau}
m&m&m&m&b&m&m&m\\
m&m&m&m&b&m&m\\
a&a&a&a\\
m&m&m\\
m&m&m\\
m
\end{ytableau}} \end{varwidth}};
\draw[very thick, black] ([xshift=0.4em, yshift=+1.5em] n.west)--([xshift=6.5em,yshift=1.5em]n.west)--([xshift=0em,yshift=-.4em] n.north);

\end{tikzpicture}
\endpgfgraphicnamed

\end{minipage}

\end{center}

\end{example}

\ytableausetup{boxsize=1.75em}
It is sometimes useful to fuse the parts of $\rho$ and $\sigma$ together as parts of size at least $m*\nu(\gamma)$ and congruent to $b\pmod m$ by appending each part of $\sigma$ at the end of each part of $\rho$. We use $\rho|\sigma$ to denote the partition obtained by this fusion and call it the \emph{enlarged sky}. To separate the partition $\rho|\sigma$ back to separate partitions $\rho$ and $\sigma$, we must know the number of parts of $\gamma$.

\begin{example}  If $\rho=\{4m,4m,4m\}$ and $\sigma=\{2m+b,2m+b,b\}$, then $\rho|\sigma=\{6m+b,6m+b,4m+b\}$\\[0.5em]

\centering
\begin{ytableau}
m&m&m&m&\none &b&m&m\\
m&m&m&m&\none &b&m&m\\
m&m&m&m& \none &b\\
\end{ytableau}
$\;\xrightarrow{\hspace{0.3in}}\;$
\begin{ytableau}
b&m&m&m&m&m&m\\
b&m&m&m&m&m&m\\
b&m&m&m&m\\
\end{ytableau}
\end{example}

\begin{example}
If $\rho|\sigma=\{5m+b,4m+b,4m+b,4m+b\}$ and there are $3$ ground parts, we know that $\rho=\{3m,3m,3m,3m\}$ and $\sigma=\{2m+b,m+b,m+b,m+b\}$. 
\end{example}

We now show that $\cp_{a,b,m}(n)$ is a generalization of $\mathcal{EO}^*(n)$
\begin{theorem}
$$\cp_{1,1,2}(n)=\mathcal{EO}^*(2n)$$
\end{theorem}
\begin{proof}
We provide a bijective proof here.  This theorem may also be proven by applying $q$-series techniques to the generating function provided in Section \ref{subsec:genFunc}.

Consider a $(1,1,2)$-copartition $(\gamma,\rho,\sigma)$ of size $n$. 
Duplicate each part of $\rho|\sigma$ to obtain a partition into an even number of odd parts of size $\ge 2\nu(\gamma)$.  Let $\gamma'$ be the conjugate partition of $\gamma$. Double the size of each part in $\gamma'$ to obtain a partition into even parts of size $\le 2\nu(\gamma)$ where only the largest part appears an odd number of times.  By combining the two resulting partitions, we have a partition of size $2n$ where all even parts are smaller than all odd parts and the largest even part is the only part appearing an odd number of times.

The inverse map is defined by taking a partition counted by $\mathcal{EO}^*(2n)$, taking the conjugate of the partition consisting of all the even parts divided by two as the ground and half of the appearances of each odd part as the enlarged sky. 
\end{proof}

\subsection{Copartition Generating Function}\label{subsec:genFunc}
Next, we explore the generating function for copartitions, showing that it can be written as an infinite product. 
\begin{theorem}\label{thm:sumProd}
Define $\cp_{a,b,m}(w,s,n)$ to be the number of $(a,b,m)$-copartitions of size $n$ that have $w$ ground parts and $s$ sky parts. Then, 
\begin{align*}
   {\mathbf{cp}}_{a,b,m}(x,y,q)&:= \sum_{n=0}^\infty\sum_{w=0}^\infty\sum_{s=0}^\infty \cp_{a,b,m}(w,s,n)x^{s}y^{w}q^n\\
   &=\frac{(xyq^{a+b};q^m)_\infty}{(xq^b;q^m)_\infty(yq^a;q^m)_\infty}.
\end{align*}
\end{theorem}
\begin{proof}[Analytic proof] 
We begin by noting that $$\frac{q^{msw+aw+bs}}{(q^m;q^m)_{w}(q^m;q^m)_{s}}$$ generates all copartitions with $s$ sky parts and $w$ ground parts. Then, by summing over $w$ and $s$ and using the variables $x$ and $y$ to keep track of the number of sky and ground parts, respectively, we can see that $$\mathbf{cp}_{a,b,m}(x,y,q)=\sum_{w=0}^\infty\sum_{s=0}^\infty \frac{x^{s}y^{w}q^{msw+aw+bs}}{(q^m;q^m)_{w}(q^m;q^m)_{s}}.$$
Then, 
\begin{align}
\mathbf{cp}_{a,b,m}(x,y,q)&=\sum_{w=0}^\infty\sum_{s=0}^\infty \frac{x^{s}y^{w}q^{msw+aw+bs}}{(q^m;q^m)_{w}(q^m;q^m)_{s}} \nonumber\\
&=\sum_{w=0}^\infty \frac{y^wq^{aw}}{(q^m;q^m)_w}\sum_{s=0}^\infty\frac{x^sq^{msw+bs}}{(q^m;q^m)_s} \nonumber\\
&=\sum_{w=0}^\infty \frac{y^{w}q^{aw}}{(q^m;q^m)_{w}(xq^{mw+b};q^m)_\infty} \label{qbinom1}\\
&=\frac{1}{(xq^b;q^m)_\infty}\sum_{w=0}^\infty\frac{y^wq^{aw}(xq^b;q^m)_w}{(q^m;q^m)_w} \nonumber\\
&=\frac{(xyq^{a+b};q^m)_\infty}{(xq^b;q^m)_\infty(yq^a;q^m)_\infty}\label{qbinom2}.
\end{align}
Note that \eqref{qbinom1} and \eqref{qbinom2} follow from the $q$-binomial theorem \cite[Eq.~17.2.37]{NIST:DLMF}. 
\end{proof}

\begin{proof}[Combinatorial proof of Theorem \ref{thm:sumProd}.]
Consider the following equivalent form of Theorem \ref{thm:sumProd}. 
\begin{equation}\label{eqn:adjthm}
   \frac{1}{(xq^b;q^m)_\infty(yq^a;q^m)_\infty}=\frac{1}{(xyq^{a+b};q^m)_\infty}
    \sum_{n=0}^\infty\sum_{w=0}^\infty\sum_{s=0}^\infty \cp_{a,b,m}(w,s,n)x^{s}y^{w}q^n.
\end{equation}

Let $\mathcal{P}_{r,m}$ be the set of partitions into parts congruent to $r\pmod{m}$ and let $\CP_{a,b,m}$ be the set of $(a,b,m)$-copartitions. To prove (\ref{eqn:adjthm}), we define a size-preserving bijection $$\phi:\mathcal{P}_{a,m}\times \mathcal{P}_{b,m}\to \mathcal{P}_{a+b,m}\times \CP_{a,b,m}.$$ 

Let $(\pi,\lambda)\in \mathcal{P}_{a,m}\times\mathcal{P}_{b,m}$. The main idea of the bijection is to pair parts of $\lambda$ with parts of $\pi$ until the remaining parts of $\lambda$ are all of size at least $m$ times the number of parts remaining in $\pi$. To achieve that goal, we let $k$ be the smallest integer such that
\begin{equation}\label{eq:bijthreshhold}
    (\lambda_k-b)/m+\nu(\lambda)-k<\nu(\pi).
\end{equation} Then, all the parts $\lambda_j$ of $\lambda$ with $j\ge k$ are matched with a part of $\pi$ to form a partition in $\mathcal{P}_{a+b,m}$ and the remaining parts of $\pi$ and $\lambda$ form a copartition. The details lie in our choice of which part of $\pi$ to match with each of the final $\nu(\lambda)-k+1$ parts of $\lambda$.  Specifically, for each $\lambda_j$ with $k\le j\le \nu(\lambda)$, combine $\lambda_j$ with $\pi_{\nu(\pi)-(\nu(\lambda)-j)-(\lambda_j-b)/m}$ to create a part of size $a+b\pmod{m}$. This creates a new partition in $\mathcal{P}_{a+b,m}$. Next, note that there are $\nu(\pi)-(\nu(\lambda)-k+1)$ parts of $\pi$ remaining. Moreover, by (\ref{eq:bijthreshhold}), the smallest remaining part of $\lambda$ is larger than $m(\nu(\pi)-\nu(\lambda)+k-1).$ Therefore, we are left with a copartition $(\gamma,\rho,\sigma) \in \mathcal{CP}_{a,b,m}$ where the remaining parts of $\lambda$ form $\rho|\sigma$ (the enlarged sky) and the remaining parts of $\pi$ become $\gamma$ (the ground). 

To see that $\phi$ is a bijection, we give the inverse map: $$\phi^{-1}: \mathcal{P}_{a+b,m}\times\CP_{a,b,m}\to \mathcal{P}_{a,m}\times\mathcal{P}_{b,m}.$$  Let $\mu\in \mathcal{P}_{a+b,m}$ and $(\gamma,\rho,\sigma)\in \CP_{a,b,m}$.
For each $0\le k<\nu(\mu)$, choose the smallest $j$ such that 
\begin{equation}\label{eq:invthresh}
\mu_{\nu(\mu)-k}-m(j)-b \le\gamma_{\nu(\gamma)-j}
\end{equation} and call it $j_k$ (when necessary, we define $\gamma_0=\infty$). Then, for each $0\le k<\nu(\mu)$, we add a part of size $m(j_k)+b$ to $\rho|\sigma$ and a part of size $\mu_{\nu(\mu)-k}-(m(j_k)+b)$ to $\gamma$. Then, we have $(\gamma,\rho|\sigma)\in \mathcal{P}_{a,m}\times\mathcal{P}_{b,m}$.
\end{proof}

\begin{example}
We work through an example of the bijection, starting with the pair $(\pi,\lambda)=(\{9,5^4,1^3\},\{26^3,22,6^2,2\})\in \mathcal{P}_{1,4}\times\mathcal{P}_{2,4}$. Note that $\nu(\pi)=8$ and $\nu(\lambda)=7$. First, we find the smallest $k$ satisfying (\ref{eq:bijthreshhold}). Note that $$(\lambda_4-2)/4+\nu(\lambda)-4=20/4+3\ge 8 \; \text{and}$$ $$(\lambda_5-2)/4+\nu(\lambda)-5=4/4+2<8,$$ so $5$ is the smallest $k$ satisfying (\ref{eq:bijthreshhold}).

Then, for each $5\le j\le 7$, we need to find the correct part of $\pi$ to match with $\lambda_j$:
\[\begin{array}{c|c|c|c}
j & \lambda_j&\nu(\pi)-(\nu(\lambda)-j)-(\lambda_j-b)/m  & \pi_{\nu(\pi)-(\nu(\lambda)-j)-(\lambda_j-b)/m} \\ \hline
5&6&5&5\\
6&6&6&1\\
7&2&8&1
\end{array}\]
Then, by combining the parts in the second and fourth columns, we create $\{11,7,3\}\in \mathcal{P}_{3,4}$. The remaining parts become a copartition with $\gamma=\{9,5^3,1\}$ and $\rho|\sigma = \{26^3,22\}$, so $$\phi((\{9,5^4,1^3\},\{26^3,22,6^2,2\}))=(\{11,7,3	\},(\{9,5^3,1\},\{20^4\},\{6^3,2\})).$$

For the inverse, we start with $(\mu,(\gamma,\rho,\sigma))=(\{11,7,3\},(\{9,5^3,1\},\{20^4\},\{6^3,2\}))\in \mathcal{P}_{3,4}\times\CP_{1,2,4}$. Then, for each $0\le k< 3$, we must figure out the smallest $j$ satisfying \eqref{eq:invthresh}:
\[ \begin{array}{c|c|c}
k&\mu_{\nu(\mu)-k}& j_k \\ \hline
0 & 3 & 0\\
1 & 7 & 1\\
2 & 11 & 1
\end{array}\]
Then, we add parts of size $4\cdot 0+2={\bf 2}$, $4\cdot1+2={\bf 6}$, and $4\cdot 1+2={\bf 6}$ to $\rho|\sigma$ and parts of size $3-2={\bf 1}$, $7-6={\bf 1}$, and $11-6={\bf 5}$ to $\gamma$. We obtain $\phi^{-1}((\{11,7,3\},(\{9,5^3,1\},\{20^4\},\{6^3,2\})))=(\{9,5^3,{\bf 5},1,{\bf 1}^2\},\{26^3,22,{\bf 6}^2,{\bf 2}\})$, where the bold parts are those that were created from $\mu=\{11,7,3\}$. 

\end{example}

\begin{example} The graphical representations of the partitions provide an insightful interpretation of the bijection $\phi$ from the proof of Theorem \ref{thm:sumProd}. We show this graphical interpretation through an example. For this example, we start with the pair $$(\{4m+a,3m+a,3m+a,2m+a,a\},\{5m+b,4m+b,2m+b,2m+b\})\in \mathcal{P}_{a,m}\times \mathcal{P}_{b,m}.$$ \\[1em]

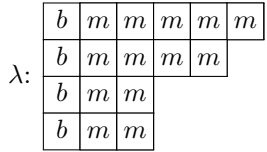
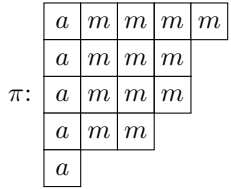
\captionsetup{type=figure}
\ytableausetup{boxsize=1.35em}
\subcaptionbox{Draw $\lambda$ and $\pi$ as $m$-modular diagrams. \label{subfig:lambdapi}}[0.47\textwidth][l]{ $\lambda$: \begin{minipage}{0.25\textwidth}
\ydiagram[*(white)m]{1+5,1+4,1+2,1+2}*[*(white)b]{6,5,3,3}\end{minipage}\\[1em]
\noindent $\pi$: \begin{minipage}{0.25\textwidth}
\ydiagram[*(white)m]{1+4,1+3,1+3,1+2,0}*[*(white)a]{1,1,1,1,1}
\end{minipage}}
\hfill
\subcaptionbox{Rotate $\pi$ clockwise by $90^\circ$ and shift each $\lambda_i$  $\nu(\lambda)-i$ squares to the right. This skewing provides a graphical interpretation of the left side of \eqref{eq:bijthreshhold}.\label{subfig:altered}}[0.47\textwidth]{ \begin{minipage}{0.40\textwidth}
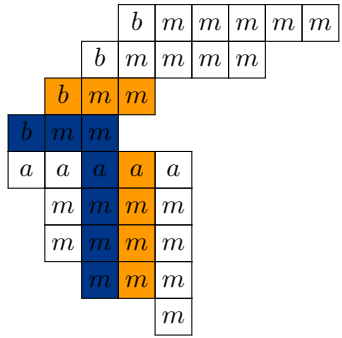
\begin{ytableau}
\none&\none&\none&b&m&m&m&m&m\\
\none&\none&b&m&m&m&m\\
\none&b&m&m\\
b&m&m\\
\end{ytableau}
\end{minipage} \\[1em]
\begin{minipage}{0.40\textwidth}
\begin{ytableau}
a&a&a&a&a\\
\none&m&m&m&m\\
\none&m&m&m&m\\
\none&\none&m&m&m\\
\none&\none&\none&\none&m\\
\end{ytableau}
\end{minipage}}
\subcaptionbox{Combine each part of $\lambda$ with the part of $\pi$ lying below its final cell, if such a part exists.\label{subfig:glued}}[0.47\textwidth]{\begin{ytableau}
\none&\none&\none&b&m&m&m&m&m\\
\none&\none&b&m&m&m&m\\
\none&*(lorange)b&*(lorange)m&*(lorange)m\\
*(dblue)b&*(dblue)m&*(dblue)m\\
a&a&*(dblue)a&*(lorange)a&a\\
\none&m&*(dblue)m&*(lorange)m&m\\
\none&m&*(dblue)m&*(lorange)m&m\\
\none&\none&*(dblue)m&*(lorange)m&m\\
\none&\none&\none&\none&m\\
\end{ytableau}}\hfill
\subcaptionbox{Separate into a partition in $\mathcal{P}_{a+b,m}$ and a copartition in $\CP_{a,b,m}$.\label{subfig:gluedsep}}[0.47\textwidth]{  \begin{minipage}{0.35\textwidth}
\vspace{5pt} \begin{ytableau}
  *(lorange)a&*(lorange)b&*(lorange)m&*(lorange)m&*(lorange)m&*(lorange)m&*(lorange)m\\
 *(dblue)a&*(dblue)b& *(dblue)m&*(dblue)m&*(dblue)m&*(dblue)m&*(dblue)m\\
 \end{ytableau}
\end{minipage}\\[1em]
\begin{minipage}{0.35\textwidth}
\begin{ytableau}
m&m&m&b&m&m\\
m&m&m&b&m\\
a&a&a\\
m&m\\
m&m\\
m\\
m\\
\end{ytableau}\end{minipage}}

\end{example}

\ytableausetup{boxsize=1.3em}

\begin{remark}
From Theorem \ref{thm:sumProd}, we can see that $(a,b,m)$-copartitions are equinumerous with the $(m,a,b)$-capsid partitions defined by Garvan and Schlosser in \cite{Capsids}.
\end{remark}

\begin{remark}
From both the graphical representations of copartitions and the generating function, we can see that $$\cp_{sa,sb,sm}(sn)=\cp_{a,b,m}(n)$$ for any $s,a,b,m\in \N$ and $n\in \N_0$. Thus, we focus our study of $\cp_{a,b,m}(n)$ on cases where $\gcd(a,b,m)=1$.
\end{remark}

By Theorem \ref{thm:sumProd}, we can see that $\mathbf{cp}_{a,b,m}(1,1,q)=\mathbf{cp}_{b,a,m}(1,1,q)$.  One reason copartitions are a logical combinatorial interpretation of the coefficients of $$\dfrac{(q^{a+b};q^m)_\infty}{(q^a;q^m)_\infty(q^b;q^m)_\infty}$$ is that this symmetry presents itself clearly through conjugation.  Recall that all partitions $\lambda$ have a conjugate partition $\lambda'$, which is obtained by reflecting the Young diagram of $\lambda$ about the line $y=-x.$ 
Similarly, we define the conjugate of a copartition $(\gamma,\rho,\sigma)$ as the copartition obtained by reflecting the graphical representation about the line $y=-x$. This conjugate copartition is precisely $(\sigma,\rho',\gamma),$ where $\rho'$ consists of exactly $\nu(\gamma)$ parts of size $m\nu(\sigma)=m\nu(\rho)$. Equivalently, $\rho'$ is the partition obtained by conjugating the $m$-modular diagram representing $\rho$.

\begin{figure}[htb]
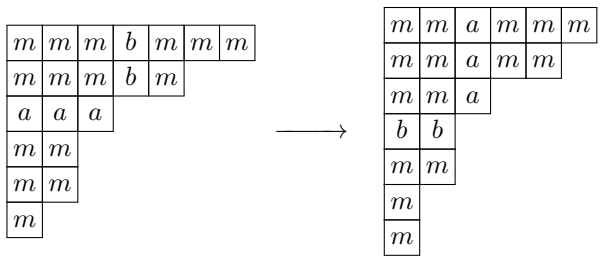

\centering
\begin{minipage}{10em}
      \begin{ytableau}
        m&m&m&b&m&m&m\\
        m&m&m&b&m\\
        a&a&a\\
        m&m\\
        m&m\\
        m
        \end{ytableau}  
\end{minipage}
   $\hspace{-2pt}\xrightarrow{\hspace{0.3in}
   }\quad$
   \begin{minipage}{10em}
        \begin{ytableau}
        m&m&a&m&m&m\\
        m&m&a&m&m\\
        m&m&a\\
        b&b\\
        m&m\\
        m\\
        m
        \end{ytableau}
        \end{minipage}
    \caption{Conjugation of an $(a,b,m)$-copartition.}
    \label{fig:my_label}
\end{figure}

\begin{remark}
Conjugation is a size-preserving bijection from $\CP_{a,b,m}$ to $\CP_{b,a,m}$.
\end{remark}

\begin{remark}
In \cite{Capsids}, Garvan and Schlosser introduce an equivalent conjugation operation on capsid partitions, noting that the symmetry is ``not at all combinatorially obvious." 
In contrast, the symmetry of copartitions is clear.

\end{remark}

\section{Special Cases}\label{sec:cop:cases}

In this section, we explore special cases of $(a,b,m)$-copartitions, highlighting connections with the Rogers-Ramanujan Functions, the sum of the partition function, and theta functions.
\subsection{Rogers-Ramanujan Functions}
Two well-studied functions in the areas of partitions and $q$-series are the Rogers-Ramanujan functions \cite[Sec.~17.2(vi)]{NIST:DLMF}: 
\begin{align}
G(q):=\sum_{n=0}^\infty\frac{q^{n^2}}{(q;q)_n}&=\frac{1}{(q;q^5)_\infty(q^4;q^5)_\infty}\label{eq:RR1}\\
H(q):=\sum_{n=0}^\infty \frac{q^{n^2+n}}{(q;q)_n}&=\frac{1}{(q^2;q^5)_\infty(q^3;q^5)_\infty}\label{eq:RR2}.
\end{align}
These functions were originally studied by Rogers \cite{Rogers}, independently discovered by Ramanujan and Schur \cite{Schur}, and further studied in a paper by both Rogers and Ramanujan \cite{RogersRamanujan}. Each of these functions have both a sum and a product representation.  The sum representations highlight the connections to basic hypergeometric series, while the product representations highlight the connections to modular functions.

Note that ${q^{n^2}}/{(q;q)_n}$ is the generating function for partitions into exactly $n$ parts such that adjacent parts have difference at least 2 (otherwise known as $2$-distinct partitions),
and $\frac{1}{(q;q^5)_\infty(q^4;q^5)_\infty}$ is the generating function for partitions into parts of size $1,4\pmod5$. Therefore, we can rephrase (\ref{eq:RR1}) as the following theorem. 
\begin{theorem}[Rogers-Ramanujan]
The number of $2$-distinct partitions of $n$ is equal to the number of partitions of $n$ into parts that are congruent to either $1$ or $4$ modulo $5$. 
\end{theorem}

Similarly, we can combinatorially state (\ref{eq:RR2}) as follows.

\begin{theorem}[Rogers-Ramanujan]
The number of $2$-distinct partitions of $n$ with no parts of size $1$ is the same as the number of partitions of $n$ into parts that are congruent to either $2$ or $3$ modulo $5$. 
\end{theorem}

Through Theorem \ref{thm:sumProd} and  the product sides of (\ref{eq:RR1}) and (\ref{eq:RR2}), we can see the following connection between the Rogers-Ramanujan functions and copartitions.
\begin{cor}\label{cor:RRConnection}
If $G(q)$ and $H(q)$ are the Rogers-Ramanujan functions defined in \eqref{eq:RR1} and \eqref{eq:RR2}, then
\begin{align*}
G(q)&=\frac{1}{(q^5;q^5)_\infty}\sum_{n=0}^\infty\cp_{1,4,5}(n)q^n \quad \quad \text{and}\\
H(q)&=\frac{1}{(q^5;q^5)_\infty}\sum_{n=0}^\infty\cp_{2,3,5}(n)q^n.
\end{align*}
\end{cor}

Corollary \ref{cor:RRConnection} motivates directions for future research. For example, a bijective proof of Corollary \ref{cor:RRConnection} could illuminate Rogers-Ramanujan function combinatorics.


\subsection{$(1,1,1)$-copartitions}

Another fascinating special case 
comes from setting $a=b=m=1$.
Since the diagrams of $(1,1,1)$-copartitions have a $1$ in every cell, they look like the Young diagrams of ordinary partitions.
When $a, b,$ and $m$ are not all equal, the diagram of an $(a,b,m)$-copartition uniquely determines the copartition;
however, when $a=b=m$, this is not the case.
Since every $(1,1,1)$-copartition has a specified number of ground parts and sky parts, 
the Young diagram of an ordinary partition will appear as the diagram of several different $(1,1,1)$-copartitions.

\begin{example}\label{ex:noBij}
The Young diagram of  
$5+5+4+4+4+2+2$
appears as four different (1,1,1)-copartitions.

\begin{minipage}{1 \textwidth}
\begin{center}
\begin{minipage}{9em}
\beginpgfgraphicnamed{copartitions-f2}
\begin{tikzpicture}
\node (n) {$\ydiagram[*(white)1]{5,5,4,4,4,2,2}$};
\draw[very thick, black] ([xshift=.4em, yshift=-4.7em] n.west)--([xshift=-3.4em,yshift=-.4em] n.north);
\end{tikzpicture}
\endpgfgraphicnamed
\end{minipage}
\hspace{1em}
\begin{minipage}{9em}
\beginpgfgraphicnamed{copartitions-f3}
\begin{tikzpicture}
\node (m) {$\ydiagram[*(white)1]{5,5,4,4,4,2,2}$};
\draw[very thick, black] ([xshift=0.4em, yshift=-2em] m.west)--([xshift=3em,yshift=-2em]m.west)--([xshift=-.7em,yshift=-.4em] m.north);
\end{tikzpicture}
\endpgfgraphicnamed
\end{minipage}
\hspace{1em}
\begin{minipage}{9em}
\beginpgfgraphicnamed{copartitions-f4}
\begin{tikzpicture}
\node (p) {$\ydiagram[*(white)1]{5,5,4,4,4,2,2}$};
\draw[very thick, black] ([xshift=0.4em, yshift=2.1em] p.west)--([xshift=5.7em,yshift=2.1em]p.west)--([xshift=2em,yshift=-.4em] p.north);
\end{tikzpicture}
\endpgfgraphicnamed
\end{minipage}
\hspace{1em}
\begin{minipage}{9em}
\beginpgfgraphicnamed{copartitions-f5}
\begin{tikzpicture}
\node (q) {$\ydiagram[*(white)1]{5,5,4,4,4,2,2}$};
\draw[very thick, black] ([xshift=0.4em, yshift=4.7em] q.west)--([xshift=7.1em,yshift=4.7em]q.west);
\end{tikzpicture}
\endpgfgraphicnamed
\end{minipage}

$(\emptyset,\emptyset,5^2+4^3+2^2)$ \qquad \qquad
$(2^2,2^5,3^2+2^3)$ \qquad  \qquad $(5^2+3^2,4^2,1^2)$ \qquad \qquad  \ $(7^2+5^2+2,\emptyset,\emptyset)$\end{center}
\end{minipage}

\end{example}

\begin{remark}
Example \ref{ex:noBij} demonstrates the more general truth that the Young diagram of an ordinary partition $\lambda$ can be realized as the diagram of a $(1,1,1)$-copartition in $\mathrm{dv}(\lambda) + 1$ ways.
\end{remark}

At the same time, it is not hard to see that the Young diagram along with the number of ground parts does completely determine the corresponding $(1,1,1)$-copartition. 
This leads us to the following.

\begin{theorem}\label{thm:cp111formula}
For all $n\in \N_0$, 
$$\cp_{1,1,1}(n)=\sum_{k=0}^n p(k).$$
\end{theorem}
\begin{proof}
We prove a bijection between the partitions of $n-k$ and the $(1,1,1)$-copartitions of $n$ with $k$ ground parts. 

Let $\lambda$ be a partition of size $n-k$. Choose $1\le j\le \nu(\lambda)$ such that $\lambda_j$ is the largest part that is no larger than $k$ or choose $j=\nu(\lambda)+1$ if all parts of $\lambda$ are larger than $k$. Then, we can define a $(1,1,1)$-copartition $(\gamma,\rho,\sigma)$ such that $\rho|\sigma=\{\lambda_1,\lambda_2,\ldots,\lambda_{j-1}\}$ and $\gamma'=\{k,\lambda_j,\lambda_{j+1},\ldots,\lambda_{\nu(\lambda)}\}$. 

To see that this map is a bijection, we note that the inverse map is defined by taking the graphical representation of a $(1,1,1)$-copartition, subtracting $1$ from each ground part, and reading the remainder of the diagram as an ordinary Young diagram for a partition of size $n-k$. 
\end{proof}

One can also give a simple proof via generating functions using Theorem \ref{thm:sumProd}.
However, the proof we give above demonstrates that the truth of Theorem \ref{thm:cp111formula} is much more fundamental, as it relies only on the definitions.

\begin{example}
For the partition $\{8,6,5,3\}$ of $22=27-5$, we have $j=3$. 
Thus $\rho | \sigma = 8+6$ and $\gamma' = 5+5+3$, and so we obtain the associated copartition $(3^3+2^2,5^2,3+1)$. 
\end{example}

The partition sum in Theorem \ref{thm:cp111formula} appears in the Online Encyclopedia of Integer Sequences \cite{OEISpartialSums} as sequence A000070. From the equivalent sequences listed there, we obtain the following corollaries. 
\begin{cor}\label{cor:partsofsize1}
The number of $(1,1,1)$-copartitions of $n-1$
is exactly the number of parts of size $1$ among all ordinary partitions of $n$.
\end{cor}

\begin{cor}\label{cor:diversities}
The number of $(1,1,1)$-copartitions of $n-1$
is exactly sum of the diversities of all ordinary partitions of $n$.
\end{cor}

\begin{remark}
To put the size of $\cp_{1,1,1}(n)$ in context, Theorem \ref{thm:cp111formula} and Corollary \ref{cor:partsofsize1} imply that $p(n) \le \cp_{1,1,1}(n) \le \mathrm{spt}(n+1)$, where $\mathrm{spt}(n)$ is the function that counts the number of appearances of the smallest part in all partitions of $n$ \cite{spt}. 
\end{remark}
One can prove both Corollary \ref{cor:partsofsize1} and \ref{cor:diversities} with direct bijections, although we omit those details here.


\subsection{$(0,1,1)$-copartitions}

In the definition of a copartition, if we allow parts of size zero in the ground or sky by letting $a$ or $b$ be zero, there may be infinitely many copartitions of a fixed $n$.
For example, if we set $a=0$, when $\sigma = \emptyset$, for every ordinary partition $\gamma$ of $n$
the copartition $(\widehat\gamma,\emptyset,\emptyset)$ is a copartition of $n$ for every $\widehat\gamma$ equal to $\gamma$ with arbitrarily many parts of size zero added.

Instead, we define
a $(0,b,m)$-copartition to be a triple of partitions $(\gamma, \rho, \sigma)$ as before, but we add the condition that $\sigma \not = \emptyset$ so that we may allow $\gamma$ to have parts of size $0$.
In constructing the the diagrams of $(0,b,m)$-copartitions, we include a row of zeros in the same way that we include a row of $a$s when $a \not = 0$.

\begin{example}
Below is the diagram of the $(0,1,2)$-copartition
$(2^3+0,8^4,7+3+1^2)$. 
\begin{center}
\beginpgfgraphicnamed{copartitions-f6}
\begin{tikzpicture}
\begin{minipage}{11em}
\node (n) {$\ydiagram[*(white)2]{4,4,4,4,0,0}*[*(white)1]{5,5,5,5,0,0}*[*(white)0]{5,5,5,5,4,0}*[*(white)2]{8,6,5,5,4,3}$};
\draw[very thick, black] ([xshift=0.4em, yshift=-1.4em] n.west)--([xshift=5.8em,yshift=-1.4em]n.west)--([xshift=0.01em,yshift=-.4em] n.north);

\end{minipage} 
\end{tikzpicture}
\endpgfgraphicnamed
\end{center}
\end{example}

By considering the generating function, we arrive at a nice formula for $\cp_{0,1,1}(n)$.

\begin{theorem}\label{thm:cp011formula}
For all $n\in \N_0$, 
\begin{equation}
    \cp_{0,1,1}(n)=\sum_{k=0}^{n-1} p(k)d(n-k),
    \label{eq:pd}
\end{equation}
where $d(n)$ is the number of divisors of $n$.
\end{theorem}

\begin{proof}
Redefining $\cp_{0,1,1}$ to match the definition above (nonempty sky), evaluating the generating function at $x=y=1$, and applying the $q$-binomial theorem twice, we have

\begin{align}
\mathbf{cp}_{0,1,1}(1,1,q)&=\sum_{s=1}^\infty\sum_{w=0}^\infty \frac{q^{sw+s}}{(q;q)_{w}(q;q)_{s}} =\sum_{s=1}^\infty \frac{q^{s}}{(q;q)_s}\sum_{w=0}^\infty\frac{q^{sw}}{(q;q)_w} =\sum_{s=1}^\infty \frac{q^{s}}{(q;q)_{s}(q^{s};q)_\infty} \nonumber\\
&=\frac{1}{(q;q)_\infty}\sum_{s=1}^\infty\frac{q^s(q;q)_{s-1}}{(q;q)_s} =\frac{1}{(q;q)_\infty}\sum_{s=1}^\infty\frac{q^s}{(1-q^s)} =\sum_{n=0}^\infty
p(n) q^n  \sum_{m=1}^\infty
d(m) q^m. \nonumber
\end{align}

\end{proof}

It is known that the right-hand side of \eqref{eq:pd} is equal to the total number of parts among all partitions of $n$, which is also equal to the sum of the largest parts among all partitions of $n$ \cite{oeisA006128}.
Traditionally, this fact has been approached using generating functions \cite{knopfmacher2005identities}.  
However, the combinatorics behind this is enlightening and intimately related to the associated copartitions, and so we now treat the relationship between $\cp_{0,1,1}(n)$ and the sum of the largest parts among all partitions of $n$ directly.

Connecting the diagrams of $(0,1,1)$-copartitions to the Young diagrams of ordinary partitions is again revealing, as it was in the last subsection.
If we disregard the cells containing $0$s, the remaining diagrams of $(0,1,1)$-copartitions have a $1$ in every cell, and thus they look like the Young diagrams of ordinary partitions.
However, these remaining Young diagrams do not uniquely determine the $(0,1,1)$-copartition. 
Instead, we find that
the Young diagram of each nonempty ordinary partition will appear as the diagram of several different $(0,1,1)$-copartitions.

\begin{example}\label{ex:noBij2}
All three of $(0,1^6,7+5+4+4+2+2)$, $(2^3,3^4,5+3+2^2)$, and $(2^3+0,4^4,4+2+1^2)$ yield the same diagram when the cells containing $0$s are removed. Below, we show those diagrams with $\rho$ outlined for each of the three copartitions. 

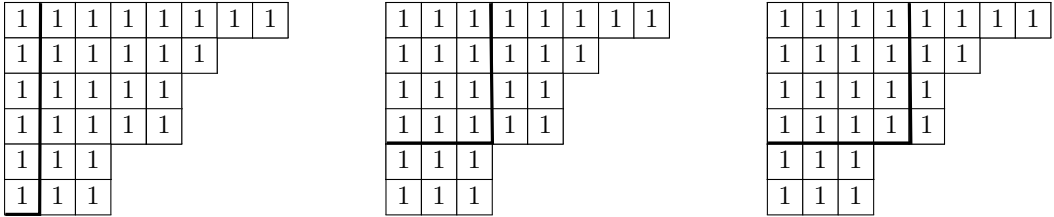
\begin{center}
\begin{minipage}{11em}
\beginpgfgraphicnamed{copartitions-f7}
\begin{tikzpicture}
\node (n) {$\ydiagram[*(white)1]{8,6,5,5,3,3}*[*(white)0]{8,6,5,5,3}$};
\draw[very thick, black] ([xshift=.4em, yshift=-4em] n.west)--([xshift=1.7em,yshift=-4em]n.west)--([xshift=-4em,yshift=-.4em] n.north);
\end{tikzpicture}
\endpgfgraphicnamed
\end{minipage}
\hspace{2em}
\begin{minipage}{11em}
\beginpgfgraphicnamed{copartitions-f8}
\begin{tikzpicture}
\node (m) {$\ydiagram[*(white)1]{8,6,5,5,3,3}*[*(white)0]{8,6,5,5,3}$};
\draw[very thick, black] ([xshift=0.4em, yshift=-1.3em] m.west)--([xshift=4.4em,yshift=-1.3em]m.west)--([xshift=-1.4em,yshift=-.4em] m.north);
\end{tikzpicture}
\endpgfgraphicnamed
\end{minipage}
\hspace{2em}
\begin{minipage}{11em}
\beginpgfgraphicnamed{copartitions-f9}
\begin{tikzpicture}
\node (p) {$\ydiagram[*(white)1]{8,6,5,5,3,3}*[*(white)0]{8,6,5,5,3}$};
\draw[very thick, black] ([xshift=0.4em, yshift=-1.3em] p.west)--([xshift=5.8em,yshift=-1.3em]p.west)--([xshift=0.01em,yshift=-.4em] p.north);
\end{tikzpicture}
\endpgfgraphicnamed
\end{minipage}   
\end{center}
\end{example}

In Theorem \ref{thm:cp011meaning}, we quantify how many different $(0,1,1)$-copartitions yield each ordinary Young diagram when removing the $0$s. 

\begin{theorem}\label{thm:cp011meaning}
For all $n\in \N_0$, 
$\cp_{0,1,1}(n)$ is the total number of parts among all partitions of $n$.
Equivalently, 
$\cp_{0,1,1}(n)$ is the sum over all partitions $\lambda$ of $n$ of the largest part of $\lambda$.
\end{theorem}

\begin{proof}
of Theorem \ref{thm:cp011meaning} from Theorem \ref{thm:cp011formula}.
The equivalence of the two statements in the theorem hold by conjugation.
Here we show that 
$\sum_{k=0}^{n-1} p(k)d(n-k)$ is the sum over all partitions $\lambda$ of $n$ of the largest part of $\lambda$.

Notice that the largest part $\lambda_1$ of any partition $\lambda = \lambda_1 + \lambda_2 + \dots + \lambda_j$ is the sum of the differences of the consecutive parts of $\lambda$,
$\lambda_1 = (\lambda_1 - \lambda_2) + (\lambda_2 - \lambda_3) + \dots + (\lambda_j - 0)$, where we consider the $j$th difference of consecutive parts to be $(\lambda_j - 0) = \lambda_j$,
and thus the sum over all partitions $\lambda$ of $n$ of the largest part of $\lambda$
is equal to the sum over all consecutive part differences among all partitions of $n$.

Notice that increasing the first $w$ parts of an arbitrary unrestricted partition of $n-hw$ by $h$ produces an arbitrary partition of $n$ with $w$th difference of consecutive parts greater than or equal to $h$.
Thus partitions of $n$ that have $w$th difference of consecutive parts greater than or equal to $h$ are in bijection with unrestricted partitions of $n-hw$.
Hence if we fix $w$ and sum $p(n-hw)$ over all $h$, we count partitions of $n$ that have $w$th difference of consecutive parts \emph{exactly} equal to $m$ $m$ times; in other words, we sum all the $w$th part differences among all partitions of $n$.
Thus, summing over all $w$, we sum all consecutive part differences among all partitions of $n$.
In doing so, we sum $p(n-hw)$ exactly $d(hw)$ times, and so   
the sum over all consecutive part differences among all partitions of $n$ is $\sum_{k=0}^{n-1} p(k)d(n-k)$ as desired.
\end{proof}

One may also prove Theorem \ref{thm:cp011meaning} combinatorially without appealing to Theorem \ref{thm:cp011formula} by constructing a direct bijection between the columns of ordinary partitions and $(0,1,1)$-copartitions of $n$, although we omit the details here.
The key element of the bijection is that each column of an ordinary partition becomes a column directly to the right of a rectangle $\rho$ of a $(0,1,1)$-copartition.

We also remark that the fact that $b=m=1$ does not play an essential role
in the proof of Theorem \ref{thm:cp011formula}.
For arbitrary $b$ and $m$, the same reasoning gives the following formula for $\cp_{0,b,m}(n)$.

\begin{theorem}\label{thm:cp0bmformula}
For all $n\in \N_0$, 
$$\cp_{0,b,m}(n)=\sum_{k<n/m} p(k)d_{b,m}(n-mk),$$
where $d_{b,m}(n)$ is the number of divisors of $n$ congruent to $b \pmod m$. 
\end{theorem}


\subsection{$(0,0,1)$-copartitions}

Following the reasoning at the beginning of the previous subsection, we define
a $(0,0,m)$-copartition to be a triple of partitions $(\gamma, \rho, \sigma)$ as before, with the added conditions that $\gamma \not = \emptyset$ and $\sigma \not = \emptyset$.

Instead of looking at the generating function, in this subsection we start with a direct bijection that gives us a formula for $\cp_{0,0,1}(n)$.
 Following \cite{straub2016core}, we define the {\em perimeter} of a partition to be the size of the largest part plus the number of parts minus one, or equivalently, the number of cells in the first row and first column of the Young diagram. 

\begin{theorem}\label{thm:cp001meaning}
For all $n\in \N_0$, 
$\cp_{0,0,1}(n)$ 
is the sum of the perimeters of the partitions of $n$.
\end{theorem}

\begin{proof}
Note that the \emph{rim} of a Young diagram is the collection of cells with no cell diagonally below and to their right.
We give a direct bijection between the cells in rims of the Young diagrams of the ordinary partitions of $n$ and $(0,0,1)$-copartitions of $n$.

 For any ordinary partition $\lambda$ and any rim cell $c$ in $\lambda$, we construct the diagram of the copartition $(\gamma, \rho, \sigma)$.
 First take the Young diagram of $\lambda$ and fill each cell with a $1$.
 Then add a row of cells containing $0$s one row below $c$ and a column of cells containing $0$s one column to the right of $c$. 
 Thus $c$ becomes the lower-right corner of $\rho$,
 each $0$ in the row of $0$s becomes the beginning of a part in $\gamma$, and each $0$ in the column of $0$s becomes the beginning of a part in $\sigma$.
 We can see that for each rim cell $c$, this generates a distinct copartition, and any $(0,0,1)$-copartition can be generated uniquely in this way. 

Since the number of cells in the rim is equal to the perimeter, our result follows.
\end{proof}
\begin{example}
Here, we show the Young diagram of $(8+6+5+5+3+3)$ with a marked rim cell and the corresponding $(0,0,1)$-copartition diagram. \\

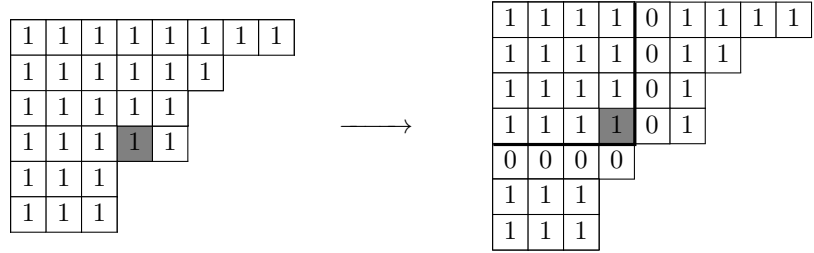
\begin{center}

\begin{minipage}{12em}
\beginpgfgraphicnamed{copartitions-f10}
\begin{tikzpicture}
\node (m) {$\ydiagram[*(white)1]{8,6,5,3,3,3}*[*(gray)1]{8,6,5,4,3,3}*[*(white)1]{8,6,5,5,3,3}$};
\end{tikzpicture}
\endpgfgraphicnamed
\end{minipage}
$\quad  \xrightarrow{\hspace{0.3in}} \quad \quad$
\begin{minipage}{13em}
\beginpgfgraphicnamed{copartitions-f11}
\begin{tikzpicture}
\node (n) {$\ydiagram[*(white)1]{4,4,4,3,0,3,3}*[*(gray)1]{4,4,4,4,0,3,3}*[*(white)0]{5,5,5,5,4,3,3}*[*(white)1]{9,7,6,6,3,3}$};
\draw[very thick, black] ([xshift=0.4em, yshift=-0.7em] n.west)--([xshift=5.75em,yshift=-0.7em]n.west)--([xshift=-0.65em,yshift=-.4em] n.north);
\end{tikzpicture}
\endpgfgraphicnamed
\end{minipage}
\end{center}
\end{example}

Since the perimeter of a partition is equal to the size of the largest part plus the number of parts minus one, Theorem \ref{thm:cp011meaning} and Theorem \ref{thm:cp001meaning} allow us to express $\cp_{0,0,1}(n)$ in terms of $\cp_{0,1,1}(n)$ and $p(n)$, which leads us to the following formulas for $\cp_{0,0,1}(n)$.

\begin{theorem}\label{thm:cp001formula}
For all $n\in \N_0$, 
$$\cp_{0,0,1}(n)= 2\cp_{0,1,1}(n) - p(n) = - p(n) + 2\sum_{k=0}^{n-1} p(k)d(n-k),$$
where $d(n)$ is the number of divisors of $n$.
\end{theorem}

\subsection{$(a,m-a,m)$-copartitions}\label{a,m-a}
When $b=m-a$, $\mathbf{cp}_{a,b,m}(x,y,q)$ is related to functions of significant analytic and combinatorial interest.  
Analytically, our interest in this case stems from the fact that we can write the generating function as a quotient of eta and theta functions.

For $|ab|<1$, one can define Ramanujan's theta function as 
$$f(a,b)=\sum_{n=-\infty}^{\infty} a^{n(n+1)/2}b^{n(n-1)/2}.$$ Note that, with an appropriate change in variables, Ramanujan's theta function is equivalent to the Jacobi theta function.
The Jacobi triple product identity tells us that $$f(a,b)=(-a;ab)_\infty(-b;ab)_\infty(ab;ab)_\infty.$$

Additionally, the Dedekind $\eta$ function is defined for $|q|<1$ as $$\eta(q)=q^{\frac{1}{24}}(q;q)_\infty.$$ The Dedekind $\eta$ function is a modular form of weight $1/2$ and level $1$. 

Using these definitions, we can see that 
\begin{align}
\mathbf{cp}_{a,m-a,m}(q)&=\frac{(q^m;q^m)_\infty}{(q^a;q^m)_\infty(q^{m-a};q^m)_\infty}\nonumber\\
&=\frac{(q^m;q^m)_\infty^2}{(q^a;q^m)_\infty(q^{m-a};q^m)_\infty(q^m;q^m)_\infty}\nonumber\\
&=q^{-\frac{m}{12}}\frac{\eta(q^m)^2}{f(-q^a,-q^{m-a})}.\nonumber
\end{align}

This connection to eta and theta functions, suggests that it should be possible to find asymptotics for the case $a+b=m$ by using similar ideas to those used in \cite{gritsenko2019theta}.

This case is also interesting combinatorially. Specifically, in \cite{Capsids}, Garvan and Schlosser used $(m,a)$-capsid partitions, which are equinumerous to $(a,m-a,m)$-copartitions in a combinatorial interpretation of Ramanujan's tau function.

\section{Conclusion}

Our definition of copartitions evolved naturally from a study partitions with parts separated by parity.
Copartitions have connections to 
a surprising number of 
other classical topics, including Roger-Ramanujan partitions, the summatory function of $p(n)$,
the total number of parts among all partitions of $n$, mock theta functions, the relationship between $\tau(n)$ and core partitions, and even crank statistics. 
Using the language of copartitions to connect Andrews' $\mathcal{EO}^*$ partitions with the capsids of Garvan and Schlosser
reveals an 
inherent symmetry
that was previously non-obvious.
Copartitions are already a three-parameter generalizations of Andrews' $\mathcal{EO}^*$ partitions, and
the
potential number of further generalizations created by changing the restrictions on partitions in the ground and sky is virtually limitless.

Many questions remain open. As a first direction, we explore divisibility properties of $\cp_{a,b,m}(n)$ in \cite{BursonEichhornParity}, but many further questions about congruences and divisibility remain open. 
In \cite{Andrews18}, Andrews defined the even-odd crank of a partition counted by $\mathcal{EO}^*(n)$ to be the largest even part minus the number of odd parts.
Strikingly, he found that the even-odd crank witnesses the congruence $
\mathcal{EO}^*(10n+8) \equiv 0 \pmod 5
$ by separating the relevant partitions into five equinumerous sets.
Translated into the language of copartitions, the even-odd crank of a copartition
becomes simply the number of parts in the ground minus the number of parts in the sky. 
A combinatorial proof that this copartition crank witnesses the congruence $\cp_{1,1,2}(5n+4) \equiv 0 \pmod 5$ is likely quite difficult, but would be illuminating.

Additionally, in \cite{BursonEichhornWeighted}, we find that the difference between $(a,b,m)$-copartitions with an even number of ground parts and $(a,b,m)$-copartitions with an odd number of ground parts is positive surprisingly often. In that same paper, we discuss an overpartition analogue for copartitions. There are many other directions that an interested reader could explore regarding weighted counts of copartitions.

The product form of the generating function for $\cp_{a,b,m}(n)$ is interesting and worthy of further study. For example, a combinatorial study of the finite products $$\frac{(q^{a+b};q^m)_n}{(q^a;q^m)_n(q^b;q^m)_n}\quad \text{and} \quad \frac{(q^{a+b};q^m)_{2n}}{(q^a;q^m)_n(q^b;q^m)_n}$$ would be compelling.

\bibliographystyle{siam}
\bibliography{copartitions}

\end{document}